\documentclass[letterpaper, 10pt, conference]{ieeeconf}

\IEEEoverridecommandlockouts
\usepackage{amsmath}
\usepackage{amssymb}

\usepackage[dvinames,table,cmyk]{xcolor}
\usepackage{tikz}
\usepackage{pgfplots}
\pgfplotsset{compat=1.18}
\usepackage{graphicx}

\usepackage{siunitx}
\usepackage{booktabs}

\usepackage[ruled,vlined,linesnumbered,noend]{algorithm2e}
\usepackage{algpseudocode}
\SetKw{KwGoTo}{go to}
\SetKw{KwTry}{try}
\SetKw{KwCatch}{catch}

\newcommand{\N}{\mathbb{N}}
\newcommand{\R}{\mathbb{R}}
\newcommand{\coloneqq}{:=}
\DeclareMathOperator*{\minimize}{minimize}
\DeclareMathOperator{\stt}{subject~to}

\newcommand{\innprod}[2]{\langle #1, #2 \rangle}

\definecolor{mydarkblue}{cmyk}{0.90,0.30,0.0,0.0}
\definecolor{mylightgreen}{cmyk}{0.30,0.0,0.45,0.0}
\definecolor{mydarkgreen}{cmyk}{0.80,0.0,1.0,0.0}

\newcommand{\GlobalPotato}{IDNP}
\newcommand{\TheTitle}{Collision Avoidance using Iterative Dynamic and Nonlinear Programming with Adaptive Grid Refinements}
\newcommand{\TheFunding}{This research work is funded by dtec.bw -- Digitalization and Technology Research Center of the Bundeswehr as part of SeRANIS -- Seamless Radio Access Networks for Internet of Space; dtec.bw is funded by the European Union -- NextGenerationEU.}

\title{\LARGE \bf \TheTitle}
\author{Rebecca Richter\and{}Alberto De~Marchi\and{}Matthias Gerdts%
\thanks{\TheFunding}%
\thanks{The authors are with the University of the Bundeswehr Munich,
	Department of Aerospace Engineering,
	Institute of Applied Mathematics and Scientific Computing,
	85577 Neubiberg, Germany
	(e-mail: \{rebecca.richter, alberto.demarchi, matthias.gerdts\}@unibw.de).
	\textcolor{red}{2024 European Control Conference (ECC), \textsc{doi}: 10.23919/ECC64448.2024.10591064.}%
    }%
}

\begin{document}

\maketitle
\thispagestyle{empty}
\pagestyle{empty}

\begin{abstract}
	Nonlinear optimal control problems for trajectory planning with obstacle avoidance present several challenges.
While general-purpose optimizers and dynamic programming methods struggle when adopted separately, their combination enabled by a penalty approach is capable of handling highly nonlinear systems while overcoming the curse of dimensionality.
Nevertheless, using dynamic programming with a fixed state space discretization limits the set of reachable solutions, hindering convergence or requiring enormous memory resources for uniformly spaced grids.
In this work we solve this issue by incorporating an adaptive refinement of the state space grid, splitting cells where needed to better capture the problem structure while requiring less discretization points overall.
Numerical results on a space manipulator demonstrate the improved robustness and efficiency of the combined method with respect to the single components.
\end{abstract}

\section{Introduction}\label{sec:intro}
Autonomous systems have been in the past decades and still are a significant subject of research in various fields, ranging from industrial applications to autonomous driving, from medical systems to space manipulators.
In this paper we consider the problem of designing collision-free trajectories for dynamical systems, possibly subject to kinematic constraints, and moving in cluttered environments.
A long-standing problem in robotics, such design is an important task in many applications \cite{hoy2015algorithms}.
Expressed in abstract terms, as a feasibility problem over the infinite-dimensional space of trajectories, the task addressed in this paper is to find states $x \colon [0,T] \to X$, controls $u \colon [0,T] \to U$, and final time $T>0$ that satisfy
\begin{align}
	&\dot{x}(t) = f(x(t),u(t)) &&\forall t \in~ [0,T] \notag\\
	&g(x(t)) \leq 0 &&\forall t \in~ [0,T] \notag \\ 
	&u(t) \in U, ~ x(t) \in X && \forall t \in~ [0,T] \label{eq:OverallProb}\\
	&b(x(T)) = 0, ~ x(0) = x_0 && \notag\\
	&T \in [T_{\min},T_{\max}]  \notag
\end{align}
with state and control bounds $X \subset \R^{n_x}$, $U \subset \R^{n_u}$, dynamics $f \colon X \times U \to \R^{n_x}$, 
nonlinear path constraints
$g \colon X \to \R^{n_g}$, final conditions $b \colon X \to \R^{n_b}$, and bounds on the final time $T_{\min}, T_{\max} > 0$.
An optimal control problem formulation is often preferred, including a cost function to filter feasible results according to application-specific preferences.
However, being obstacles the source of problematic nonconvexities in the trajectory-optimization problem, in this work we focus on reliable methods for seeking \emph{feasible} trajectories.

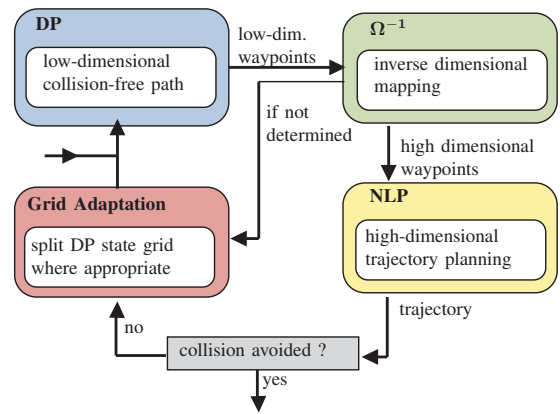
\begin{figure}[b!]
	\centering%
	\tikzset{every picture/.style={line width=0.75pt}} 

\scalebox{0.72}{%

\begin{tikzpicture}[x=0.75pt,y=0.75pt,yscale=-1,xscale=1]

\draw  [fill={rgb, 255:red, 74; green, 144; blue, 226 }  ,fill opacity=0.45 ] (99,45.4) .. controls (99,36.34) and (106.34,29) .. (115.4,29) -- (232.6,29) .. controls (241.66,29) and (249,36.34) .. (249,45.4) -- (249,94.6) .. controls (249,103.66) and (241.66,111) .. (232.6,111) -- (115.4,111) .. controls (106.34,111) and (99,103.66) .. (99,94.6) -- cycle ;
\draw [line width=1.5]    (249,70) -- (327,70) ;
\draw [shift={(331,70)}, rotate = 180] [fill={rgb, 255:red, 0; green, 0; blue, 0 }  ][line width=0.08]  [draw opacity=0] (11.61,-5.58) -- (0,0) -- (11.61,5.58) -- cycle    ;
\draw [line width=1.5]    (360,231) -- (360,273) ;
\draw [line width=1.5]    (360,273) -- (346,273) ;
\draw [shift={(342,273)}, rotate = 360] [fill={rgb, 255:red, 0; green, 0; blue, 0 }  ][line width=0.08]  [draw opacity=0] (11.61,-5.58) -- (0,0) -- (11.61,5.58) -- cycle    ;
\draw [line width=1.5]    (361,111) -- (361,146) ;
\draw [shift={(361,150)}, rotate = 270] [fill={rgb, 255:red, 0; green, 0; blue, 0 }  ][line width=0.08]  [draw opacity=0] (11.61,-5.58) -- (0,0) -- (11.61,5.58) -- cycle    ;
\draw  [fill={rgb, 255:red, 155; green, 155; blue, 155 }  ,fill opacity=0.46 ] (206,259) -- (340,259) -- (340,282) -- (206,282) -- cycle ;
\draw [line width=1.5]    (269,282) -- (269,308) ;
\draw [shift={(270,312)}, rotate = 268.49] [fill={rgb, 255:red, 0; green, 0; blue, 0 }  ][line width=0.08]  [draw opacity=0] (11.61,-5.58) -- (0,0) -- (11.61,5.58) -- cycle    ;
\draw [line width=1.5]    (170,272) -- (170,238) ;
\draw [shift={(171,234)}, rotate = 91.51] [fill={rgb, 255:red, 0; green, 0; blue, 0 }  ][line width=0.08]  [draw opacity=0] (11.61,-5.58) -- (0,0) -- (11.61,5.58) -- cycle    ;
\draw [line width=1.5]    (170,272) -- (205,272) ;
\draw  [fill={rgb, 255:red, 255; green, 255; blue, 255 }  ,fill opacity=1 ] (107,64.4) .. controls (107,59.76) and (110.76,56) .. (115.4,56) -- (231.6,56) .. controls (236.24,56) and (240,59.76) .. (240,64.4) -- (240,89.6) .. controls (240,94.24) and (236.24,98) .. (231.6,98) -- (115.4,98) .. controls (110.76,98) and (107,94.24) .. (107,89.6) -- cycle ;
\draw  [fill={rgb, 255:red, 126; green, 211; blue, 33 }  ,fill opacity=0.46 ] (329,47.4) .. controls (329,38.89) and (335.89,32) .. (344.4,32) -- (463.6,32) .. controls (472.11,32) and (479,38.89) .. (479,47.4) -- (479,93.6) .. controls (479,102.11) and (472.11,109) .. (463.6,109) -- (344.4,109) .. controls (335.89,109) and (329,102.11) .. (329,93.6) -- cycle ;
\draw  [fill={rgb, 255:red, 255; green, 255; blue, 255 }  ,fill opacity=1 ] (339,66.4) .. controls (339,61.76) and (342.76,58) .. (347.4,58) -- (463.6,58) .. controls (468.24,58) and (472,61.76) .. (472,66.4) -- (472,91.6) .. controls (472,96.24) and (468.24,100) .. (463.6,100) -- (347.4,100) .. controls (342.76,100) and (339,96.24) .. (339,91.6) -- cycle ;
\draw  [fill={rgb, 255:red, 248; green, 231; blue, 28 }  ,fill opacity=0.46 ] (329,166.4) .. controls (329,157.89) and (335.89,151) .. (344.4,151) -- (463.6,151) .. controls (472.11,151) and (479,157.89) .. (479,166.4) -- (479,212.6) .. controls (479,221.11) and (472.11,228) .. (463.6,228) -- (344.4,228) .. controls (335.89,228) and (329,221.11) .. (329,212.6) -- cycle ;
\draw  [fill={rgb, 255:red, 255; green, 255; blue, 255 }  ,fill opacity=1 ] (338,185.4) .. controls (338,180.76) and (341.76,177) .. (346.4,177) -- (462.6,177) .. controls (467.24,177) and (471,180.76) .. (471,185.4) -- (471,210.6) .. controls (471,215.24) and (467.24,219) .. (462.6,219) -- (346.4,219) .. controls (341.76,219) and (338,215.24) .. (338,210.6) -- cycle ;
\draw  [fill={rgb, 255:red, 208; green, 2; blue, 27 }  ,fill opacity=0.51 ] (99,169.4) .. controls (99,160.89) and (105.89,154) .. (114.4,154) -- (233.6,154) .. controls (242.11,154) and (249,160.89) .. (249,169.4) -- (249,215.6) .. controls (249,224.11) and (242.11,231) .. (233.6,231) -- (114.4,231) .. controls (105.89,231) and (99,224.11) .. (99,215.6) -- cycle ;
\draw  [fill={rgb, 255:red, 255; green, 255; blue, 255 }  ,fill opacity=1 ] (106,190.4) .. controls (106,185.76) and (109.76,182) .. (114.4,182) -- (230.6,182) .. controls (235.24,182) and (239,185.76) .. (239,190.4) -- (239,215.6) .. controls (239,220.24) and (235.24,224) .. (230.6,224) -- (114.4,224) .. controls (109.76,224) and (106,220.24) .. (106,215.6) -- cycle ;
\draw [line width=1.5]    (270,190) -- (254,190) ;
\draw [shift={(250,190)}, rotate = 360] [fill={rgb, 255:red, 0; green, 0; blue, 0 }  ][line width=0.08]  [draw opacity=0] (11.61,-5.58) -- (0,0) -- (11.61,5.58) -- cycle    ;
\draw [line width=1.5]    (270,80) -- (270,190) ;
\draw    (270,80) -- (330,80) ;

\draw (107.4,33) node [anchor=north west][inner sep=0.75pt, align=left] {\textbf{ DP}};
\draw (117.4,59) node [anchor=north west][inner sep=0.75pt, align=left] {low-dimensional\\ collision-free path};
\draw (349.4,61) node [anchor=north west][inner sep=0.75pt, align=left] {inverse dimensional\\mapping};
\draw (346.4,35) node [anchor=north west][inner sep=0.75pt, align=left] {$\mathbf{\Omega ^{-1}}$};
\draw (343.2,182) node [anchor=north west][inner sep=0.75pt, align=left] {high-dimensional \\trajectory planning};
\draw (346.4,154) node [anchor=north west][inner sep=0.75pt, align=left] {\textbf{NLP}};
\draw (106.6,159) node [anchor=north west][inner sep=0.75pt, align=left] {\textbf{Grid Adaptation}};
\draw (109,188) node [anchor=north west][inner sep=0.75pt, align=left] {split DP state grid\\where appropriate};
\draw (254,40) node [anchor=north west][inner sep=0.75pt, align=left] {low-dim.\\ waypoints};
\draw (368,119) node [anchor=north west][inner sep=0.75pt, align=left] {high dimensional\\waypoints};
\draw (213,262) node [anchor=north west][inner sep=0.75pt, align=left] {collision avoided ?};
\draw (367,233) node [anchor=north west][inner sep=0.75pt, align=left] {trajectory};
\draw (271,285) node [anchor=north west][inner sep=0.75pt, align=left] {yes};
\draw (174.5,248.5) node [anchor=north west][inner sep=0.75pt, align=left] {no};
\draw (274,95) node [anchor=north west, inner sep=0.75pt, align=left] {if not\\determined};

\draw [line width=1.5]    (171,155) -- (171,110) ;
\draw [shift={(171,106)}, rotate = 90] [fill={rgb, 255:red, 0; green, 0; blue, 0 }  ][line width=0.08]  [draw opacity=0] (11.61,-5.58) -- (0,0) -- (11.61,5.58) -- cycle    ;
\draw [line width=1.5]    (120,132) -- (171,132) ;
\draw [shift={(144,132)}, rotate = 180] [fill={rgb, 255:red, 0; green, 0; blue, 0 }  ][line width=0.08]  [draw opacity=0] (11.61,-5.58) -- (0,0) -- (11.61,5.58) -- cycle    ;

\end{tikzpicture}%
}%
%
	\caption{Schematic overview of the proposed iterative scheme for motion planning in constrained environments: low-dimensional collision-free waypoints from dynamic programming (DP) are mapped back to the original dimension and tracked by a trajectory planned with nonlinear programming (NLP). Possible collisions are integrated into the DP via a penalty term and the DP state space grid is adapted, until convergence to a valid solution.}%
	\label{fig:flowchart}%
\end{figure}

The popular method of dynamic programming (DP) is able to handle the constraints \eqref{eq:OverallProb} and provide a solution under mild assumptions, but suffers from the curse of dimensionality \cite{bellman1957dynamic,bertsekas2005dynamic}.
Despite multiple attempts of lowering its computational requirements, including enhanced sampling strategies \cite{cervellera2013quasi,weber2017optimized}, adaptive state space refinements \cite{richter2023adaptive, munos2002variable}, and massive parallelisation,
DP is not yet applicable in practice to high dimensional scenarios.
In contrast, sampling-based methods \cite{elbanhawi2014sampling} such as 
Rapidly-exploring Random Trees (RRT) 
\cite{lavalle1998rapidly} show good performance on larger state space dimensions.
Therefore, these methods are widespread, especially in the field of robotics, but kinematics or dynamics are hardly incorporated \cite{hoy2015algorithms}.
Approaches based on solving large-scale nonlinear programs (or systems of equations) are also limited
\cite{gerdts2023optimal}.
General-purpose optimizers like 
Ipopt \cite{waechter2006implementation} and WORHP \cite{bueskens2013esa}
can be highly sensitive to initial guesses and provide, if any, a local solution.
Since they additionally tend to struggle with many nonconvex constraints, as those for collision avoidance, several reformulations for specific obstacle structures were presented \cite{zhang2021optimization,guthrie2022differentiable}.
Numerical approaches tailored to motion planning tasks are able to improve robustness but remain sensitive to the quality of user-defined initializations \cite{zucker2013chomp, schulman2013finding}.
To overcome these difficulties, the approaches presented in \cite{natarajan2021inter, britzelmeier2023dynamic} combine low dimensional graph-search and dynamic programming, respectively, with gradient-based optimization in a higher dimensional space, exploiting the strengths of the respective methods.
In a similar spirit, ``to take the best of the two worlds'', \cite{marcucci2022motion} brings together network and convex optimization in a mixed-integer programming formulation for a special class of motion planning problems.

\paragraph*{Motivation}
In this work we focus on the framework proposed in \cite{britzelmeier2023dynamic}, here referred to as Iterative Dynamic and Nonlinear Programming (\GlobalPotato{}),
which alternates between DP and nonlinear programming (NLP) to successively improve the design relative to the original problem \eqref{eq:OverallProb}, as follows.
Exploiting a problem-specific lower-dimensional representation of \eqref{eq:OverallProb}, the trajectory planning task is decomposed into two coupled components:
(i) the generation of a sequence of low-dimensional waypoints and (ii) the dynamically-feasible motion planning through these waypoints, via DP and NLP respectively.
Seeking a collision-free trajectory, waypoints leading to infeasible motions are penalized, thus strongly coupling the two components.
Waypoints enter into the motion planning task as path constraints on the state \cite{gerdts2023optimal}, whereas the penalty reformulation in DP follows from the theory in \cite{britzelmeier2021decomposition}.
In \cite{britzelmeier2023dynamic} the state space discretization for DP is precomputed and fixed, requiring a reasonable initial sampling to cover the obstacles.
To better capture the environment with the penalty term, this is adapted through the iterations but, because of the fixed state space discretizations, cannot capture narrow passages in between obstacles on coarse discretizations, thus inhibiting the solution process.
Moreover, this is especially an issue as high-resolution sampling results in huge computational effort for DP due to the curse of dimensionality.

\paragraph*{Contribution}
In this paper, with the goal of limiting the memory and runtime requirements on fine grids and enhancing convergence on coarse ones, we incorporate an adaptive state space refinement strategy into \GlobalPotato{}, building upon ideas presented in \cite{richter2023adaptive}.
However, instead of formulating a refinement criterion as a trigger for cell refinements, we use the collisions detected in the trajectories planned by the NLP block.
We demonstrate on a space manipulator numeric example improved robustness and efficiency of the combined method compared to both the original one presented in \cite{britzelmeier2023dynamic} and full discretization.
Figure~\ref{fig:flowchart} illustrates the overall structure of the proposed iterative scheme.
Notice how each component plays an important role:
without DP, NLP struggles with nonconvex constraints;
without NLP, DP is not practical for high dimensional systems;
without adaptive refinements, DP and NLP combined (as in \GlobalPotato{}) require fine discretizations or may not converge.

\paragraph*{Outline}
We will detail the components of our approach, sketching those of \GlobalPotato{} and highlighting the proposed changes.
Following \cite{britzelmeier2023dynamic}, the low- and high-dimensional problem representations are discussed in Sections~\ref{sec:upper} and~\ref{sec:lower}, respectively.
Then, the state space discretization and adaptive splitting are considered in Section~\ref{sec:splitting}.
Finally, numerical simulations are reported in Section~\ref{sec:numerical}.

\section{Low Dimensional Level}\label{sec:upper}
The DP method is capable of handling nonconvex constraints and providing a global solution, but suffers from the curse of dimensionality \cite{bellman1957dynamic,bertsekas2005dynamic}.
In this section we briefly recall the formulation in \cite[Section II]{britzelmeier2023dynamic} to apply DP on a low dimensional representation of \eqref{eq:OverallProb}, with state and control spaces $W \subset \R^{n_w}$ and $V \subset \R^{n_v}$, respectively, associated to a suitable user-defined mapping $\Omega \colon X \to W$.
The precise definition of $\Omega$ as well as of \emph{high}- and \emph{low}-dimensional representations are problem-dependent.
As a prominent example, a robotic manipulator with 12 degrees of freedom may be described with $n_x=12$, $n_w=3$, and $\Omega$ its forward kinematics.
The \GlobalPotato{} algorithm builds upon a 
Low-Dimensional Problem \eqref{eq:DPProb}
of the form
\begin{align}
	\label{eq:DPProb}\tag{\text{LDP}}
	\minimize_{w, v, \tau_f} ~&\mathcal{M}(w(\tau_f)) + \rho \int_{0}^{\tau_f} \mathcal{P}(w(\tau)) \mathrm{d}\tau \notag \\
	\stt ~& \dot{w}(\tau) = \bar{f}(w(\tau),v(\tau)) \notag & \forall  \tau \in [0,\tau_f] \notag \\
	& v(\tau) \in V \notag & \forall  \tau \in [0,\tau_f] \notag \\
	&w(0) = \Omega(x_0) ,~ \tau_f \in [T_{\min},T_{\max}] \notag
\end{align}
with state $w \colon [0,\tau_f] \to W$, control $v \colon [0,\tau_f] \to V$, dynamics $\bar{f} \colon W \times V \to \R^{n_w}$, and a free final time $\tau_f$.
It includes the inital condition of \eqref{eq:OverallProb} using the dimensional mapping $\Omega$, and the terminal conditions are encoded into the Mayer-type term $\mathcal{M}: W \to \R_+$,
$\mathcal{M}(w)
\coloneqq
\min_{x \in X} \left\{ \| b(x) \| \,\vert\, \Omega(x) = w \right\}$.
Additionally, the penalty function $\mathcal{P} \colon W \to \R_+$,
\begin{equation} \label{eq:penalty_function}
	\mathcal{P}(w)
	\coloneqq
	\min_{x \in X} \left\{ \| \Omega(x) - w \| + \| \max\{0, g(x)\} \| \right\}
	,
\end{equation}
accounts for the collision avoidance constraints, weighted by a penalty factor $\rho > 0$.

\subsection{Dynamic Programming}
We discretize \eqref{eq:DPProb} with respect to time to obtain a finite-dimensional transcription of this optimization problem, which will then be tackled using DP.
Let us consider a fixed number of time steps $M \in \N$, and a set of time step sequences
$\mathbb{H}
\coloneqq
\left\{
h \in \R_+^M \,\vert\, \sum_{i=0}^{M-1} h_i \in [T_{\min},T_{\max}]
\right\}$,
each defining a certain discrete time grid
$\mathbb{G}_\tau$ given by time points $\tau_i = \sum_{j=0}^i h_j$, $i = 0,\ldots,M$.
The resulting time-discretized problem can be written as
\begin{align}
	\label{eq:DiscDPProb}
	\minimize_{w_h, v_h, h} ~&\mathcal{M}(w_h(\tau_M)) + \rho \sum_{i = 0}^{M-1} h_i \mathcal{P}(w_h(\tau_i)) \\
	\stt ~& w_h(\tau_{i+1}) = \bar{f}_{h_i}(w_h(\tau_i),v_h(\tau_i)) ~ \forall \tau_i,\tau_{i+1} \in \mathbb{G}_{\tau}\notag \\
	&v_h(\tau_i) \in V \qquad \qquad \qquad \qquad ~~ \forall \tau_i \in \mathbb{G}_{\tau}\notag\\
	&w_h(\tau_0) = \Omega(x_0) ,~ h\in \mathbb{H} \notag
\end{align}
with $w_h \colon \mathbb {G}_{\tau} \to W$ and $v_h \colon \mathbb{G}_{\tau} \to V$ representing the state and control approximations.
The explicit Euler approximation of the dynamics is given by 
$\bar{f}_h(w,v)
\coloneqq
w + h \bar{f}(w,v)$.
As described in details in \cite[Section 2.B]{britzelmeier2023dynamic} the method of approximate DP is applied, based on the recursion 
\begin{align}
	&\vartheta(\tau_M,w_M)
	\coloneqq
	\mathcal{M}(w_M) \label{eq:UL_approx_value_function}  \\
	&\vartheta(\tau_\ell, w_\ell)
	\coloneqq
	\min_{v_{h,\ell} \in \mathbb{G}_v, h \in \mathbb{H}}
	\Big\{
	\rho h_\ell \mathcal{P}(w_\ell) \nonumber \\	
	&\quad+\sum_{j=1}^{2^{n_w}} \gamma_j\left(\bar{f}_{h_\ell}(w_\ell,v_{h,\ell})\right) \vartheta\left(\tau_{\ell+1},\Pi_j^{\ell+1}[\bar{f}_{h_\ell}(w_\ell,v_{h,\ell})]\right)
	\Big\} \nonumber
\end{align}
defined for initial states $z_\ell \in \mathbb{G}_w^\ell$ and times $\tau_\ell \in \mathbb{G}_{\tau}$, and control grid $\mathbb{G}_{v}$.
To evalute the value function $\vartheta$ on the state grid $\mathbb{G}_w^{\ell+1}$, we employ an arbitrary interpolation, represented by $\Pi_j^{\ell+1} \colon W \to \mathbb{G}_w^{\ell+1}$ and
$\gamma_j^{\ell+1} \colon W \to \R$,
 returning respectively the $j$-th vertex of the containing cell and its associated weight.
In practice, the property $h\in\mathbb{H}$, or $\tau_f \in [T_{\min},T_{\max}]$, can be guaranteed by setting bounds on each timestep as $T_{\min}/M \leq h_i \leq T_{\max}/M$.

We stress that \eqref{eq:DiscDPProb} stems from a time-discretization of \eqref{eq:DPProb}, whereas \cite[Eq. 6]{britzelmeier2023dynamic} can be interpreted as arising from the penalization of constraint violations on the time grid $\mathbb{G}_\tau$.
Though motivated in different ways, both variants serve their purpose for the low-dimensional planning, as reported in Section~\ref{sec:numerical}.

In contrast to the original \GlobalPotato{} of \cite{britzelmeier2023dynamic}, here we define a state grid for every time point, since they might have different structures due to the adaptive refinements presented in Section~\ref{sec:splitting} below.
The convergence towards the continuous problem attested in \cite[Section 2.B]{britzelmeier2023dynamic}, on infinitely fine time, state, and control discretizations is left untouched.
Notice that in general it can be more efficient to have time-dependent, yet adaptive state grids, because they can be much coarser than in the uniformly sampled case.
Once the optimal trajectory $\{\tau_j, w_j\}_{j = 0}^{M}$ resulting from \eqref{eq:UL_approx_value_function} has been computed, it is mapped back onto the original space of \eqref{eq:OverallProb}.
If feasible, such collision-free sequence enters the NLP block, without explicit formulation of the collision constraints; cf. Figure~\ref{fig:flowchart}.

\section{Trajectory Generation}\label{sec:lower}
While the adaptive state space refinement process requires some changes to the DP block and interfaces of \GlobalPotato{}, the trajectory planning block remains as in \cite[Section III]{britzelmeier2023dynamic}.
The High-Dimensional Problem \eqref{eq:LL_Problem}
\begin{align}
	\minimize_{x, u, T} \quad &J(x, u, T) \tag{HDP} \label{eq:LL_Problem} \\
	\stt \quad & \dot{x}(t) = f(x(t),u(t)) &&\forall~t\in[0, T]\nonumber \\
	& x(t) \in X,~u(t) \in U &&\forall~t\in[0, T] \nonumber\\
	& p_j( x(T \tau_j/\tau_M) ) = 0 &&\forall~j\in\{0,\ldots,M\} \nonumber \\
	& x(0) = x_0 \nonumber \\
	& T \in [T_{\min}, T_{\max}] , \nonumber
\end{align} 
accounts for the dynamics, initial conditions, and bounds of the original problem \eqref{eq:OverallProb}.
Additionally, state and control sets $X$ and $U$ can encode general constraints, and could also be coupled, as long as simple enough to be handled by NLP solvers (e.g., convex sets described by inequalities).
While the objective function $J \colon X \times U \times [T_{\min},T_{\max}] \to \R$ can reflect some user preferences, functions $p_j \colon X \to \R$, ${j = 0,\ldots,M}$, represent the waypoint constraints.
Their explicit formulation depends on the structure of the dimensional mapping $\Omega$; the prominent example $p_j(x) \coloneqq \Omega(x) - w_j$ is used in our numerical simulations.

Note that \eqref{eq:LL_Problem} is also a free final time problem.
However, since the timestamped waypoints $\{\tau_j, w_j\}_{j = 0}^{M}$ from \eqref{eq:DiscDPProb} satisfy the time constraints by construction of $\mathbb{H}$, it would be possible to consider it fixed here, with $T = \tau_f$.
Despite the issue of dealing with a free final time optimal control problem, we prefer the formulation \eqref{eq:LL_Problem} to equip the motion planning task with an additional degree of freedom.
Nonetheless, feasibility of \eqref{eq:LL_Problem} is still not guaranteed as the waypoints could be outside the reachable set.

In order to apply off-the-shelf NLP solvers to \eqref{eq:LL_Problem}, the trajectory generation task is transcripted into a finite-dimensional problem by discretization with respect to time.
Therefore, we introduce $N \gg M$ discrete time points $\{t_i\}_{i=0}^{N}$, respective state $\{x_i\}_{i=0}^N$ and control $\{u_i\}_{i=0}^{N-1}$ approximations, and a discrete dynamics defined with respect to an arbitrary approximation scheme.
The waypoint timestamps $\{ T\tau_j/\tau_M \}_{j=0}^M$ are embedded into the discrete time points $\{t_i\}_{i=0}^{N}$ to directly encode the waypoint constraints in the time-discretized problem.
After solving a time-discretized counterpart of \eqref{eq:LL_Problem}, the solution trajectory $\{t_i,x_i\}_{i= 0}^N$ satisfies the relevant dynamics, state and control bounds, and path constraints.
Therefore, if the generated trajectory happens to avoid obstacles, then it is a valid solution to the original problem \eqref{eq:OverallProb}, up to the NLP discretization and tolerance.
Finally, notice that \eqref{eq:LL_Problem} can be addressed by any techniques able to tackle optimal control problems with state constraints, that is, not necessarily by full discretization; see \cite{gerdts2023optimal} for other approaches.

\subsection{Space Mapping}\label{sec:mapping}
In order to connect the low- and high-dimensional problems presented above, the mapping $\Omega \colon X \to W$ is used to transfer points in the original space $X$ to their low dimensional representation in $W$.
While $\Omega$ is well-defined, e.g., via forward kinematics models,
the inverse operation is often ill-posed, requiring to determine the high dimensional representation of a low dimensional trajectory.
Problem 
\begin{align}
	\label{eq:SpaceMapping}
	\minimize_{x, \sigma}\quad& \varrho_1 \| x - x_{j-1} \|^2 + \varrho_2 \innprod{1}{\sigma} \\
	\stt\quad& \Omega(x) = w_j ,\quad
	g(x) \leq \sigma \leq 0 ,\quad
	x \in X \nonumber
\end{align} 
is solved in \GlobalPotato{} to
determine a state that belongs to $\Omega^{-1}(w_j)$, if possible,
is close to a given previous state $x_{j-1} \in X$, and has a large margin before violating the collision avoidance constraints, with weighting factors $\varrho_1, \varrho_2 \geq 0$.
As pointed out in \cite[Section IV.A]{britzelmeier2023dynamic}, solving this problem recursively starting from $x_0 \in X$ reduces the issue of multiple solutions and leads to more stable high dimensional sequences $\{\tau_j,x_j\}_{j= 0}^M$.

Although obstacle avoidance constraints are relaxed and penalized, with the goal of promoting collision-free configurations, problem \eqref{eq:SpaceMapping} can be infeasible if the waypoint $w_j$ is not reachable with the forward mapping $\Omega$.
For instance, many robotic manipulators have limited reachable regions because of their mechanical construction.
If this is the case, that is, if a waypoint $w_j$ selected by DP is deemed infeasible at this stage, then the corresponding point in the state grid $\mathbb{G}_w^j$ is penalized in \eqref{eq:DPProb} via $\mathcal{P}$ to discourage its selection at the next DP call.
A similar procedure is applied if a trajectory point $x_i$ from \eqref{eq:LL_Problem} is not collision-free. 
In \cite[Alg. 1, Step 10]{britzelmeier2023dynamic}, such penalization would be distributed onto neighboring state grid points, without affecting the state grid itself.
In contrast, here an infeasible inverse mapping problem triggers a direct call to the grid adaptive refinement, as illustrated in Figure~\ref{fig:flowchart} and discussed in the following section.

\section{Adaptive Grid Refinement}\label{sec:splitting}
Deviating from the \GlobalPotato{} algorithm of \cite{britzelmeier2023dynamic}, we replace the penalty adaptation therein with adaptive refinements of the state space grid.
In this way, instead of lumping penalties together, constraints are better captured by the penalty landscape thanks to additional evaluation points, located where needed.
In fact, the state space can initially be coarsely sampled and fine obstacle structures will then be recognized, if of interest, steering convergence to feasible trajectories. 
For reference, we apply the adaptive mesh approach from \cite{richter2023adaptive} instead of the penalty approximation in \cite[Section IV.B]{britzelmeier2023dynamic}.

As pointed out in Section~\ref{sec:upper}, we consider state space grids $\{\mathbb{G}_w^j\}_{j = 0}^M$ dividing the admissible space into cells $\{\mathbb{C}_w^j\}_{j = 0}^M$, for each discrete time step of the DP algorithm.
Allowing arbitrary grid structures, we consider some functions $C_w^j \colon W \to \mathbb{C}_w^j$ and $\mathcal{S}_w^j \colon \mathbb{C}_w^j \to P(W)$, $j=0,\ldots,M$.
The former map a certain point to its corresponding cell, while the latter are splitting operations, with potential set $P$, returning the grid extended with the points resulting from the splitting of a specific cell.
For a detailed discussion we refer to \cite[Section III]{richter2023adaptive}, where the explicit definitions for the case of $n$-dimensional cube grids are also presented.

As depicted in Figure~\ref{fig:flowchart}, the grid adaptation process is called in two cases, and with two types of input.
On the one hand, the inverse mapping hands over sets $\{\tau_k,w_k\}_{k = 0}^{K}$ of $K \in (0,M]$ waypoints deemed infeasible.
On the other hand, the NLP optimized trajectory has $K \in (0,N]$ points $\{t_k,x_k\}_{k = 0}^{K}$ that are not collision-free.
In the latter case, the points are transferred to the DP space applying the forward dimensional mapping $\Omega$, at time $\tau_k = \tau_M t_k / T$ and space $w_k = \Omega(x_k)$ respectively, such that a low dimensional sequence $\{\tau_k,w_k\}_{k = 0}^{K}$ is generated.

\begin{figure}[tbh]
	\centering%
	 \tikzset{every picture/.style={line width=0.75pt}} 
 
 \scalebox{0.8}{%
 
 \begin{tikzpicture}[x=0.75pt,y=0.75pt,yscale=-1,xscale=1]
	\draw [color={rgb, 255:red, 65; green, 117; blue, 5 }  ,draw opacity=1 ][line width=3]  [dash pattern={on 3.38pt off 3.27pt}]  (31.68,39.14) .. controls (67.6,25.66) and (94.06,151.44) .. (157.69,111.57) ;	
	\draw [color={rgb, 255:red, 208; green, 2; blue, 27 }  ,draw opacity=1 ][line width=3]  [dash pattern={on 3.38pt off 3.27pt}]  (56.26,47.56) .. controls (70,65) and (92.8,105) .. (110.52,110.63) ;
	
	\filldraw [color={rgb, 255:red, 0; green, 0; blue, 0 }] (0,0) circle (2pt);
	\filldraw [color={rgb, 255:red, 0; green, 0; blue, 0 }] (0,50) circle (2pt);
	\filldraw [color={rgb, 255:red, 0; green, 0; blue, 0 }] (0,100) circle (2pt);
	\filldraw [color={rgb, 255:red, 0; green, 0; blue, 0 }] (0,150) circle (2pt);
	\filldraw [color={rgb, 255:red, 0; green, 0; blue, 0 }] (50,0) circle (2pt);
	\filldraw [color={rgb, 255:red, 0; green, 0; blue, 0 }] (50,50) circle (2pt);
	\filldraw [color={rgb, 255:red, 0; green, 0; blue, 0 }] (50,100) circle (2pt);
	\filldraw [color={rgb, 255:red, 0; green, 0; blue, 0 }] (50,150) circle (2pt);
	\filldraw [color={rgb, 255:red, 0; green, 0; blue, 0 }] (100,0) circle (2pt);
	\filldraw [color={rgb, 255:red, 0; green, 0; blue, 0 }] (100,50) circle (2pt);	
	\filldraw [color={rgb, 255:red, 0; green, 0; blue, 0 }] (100,100) circle (2pt);
	\filldraw [color={rgb, 255:red, 0; green, 0; blue, 0 }] (100,150) circle (2pt);
	\filldraw [color={rgb, 255:red, 0; green, 0; blue, 0 }] (150,0) circle (2pt);
	\filldraw [color={rgb, 255:red, 0; green, 0; blue, 0 }] (150,50) circle (2pt);
	\filldraw [color={rgb, 255:red, 0; green, 0; blue, 0 }] (150,100) circle (2pt);
	\filldraw [color={rgb, 255:red, 0; green, 0; blue, 0 }] (150,150) circle (2pt);
	
	\filldraw [color={rgb, 255:red, 37; green, 150; blue, 190 }] (75,75) circle (1.5pt);
	\filldraw [color={rgb, 255:red, 37; green, 150; blue, 190 }] (75,50) circle (1.5pt);
	\filldraw [color={rgb, 255:red, 37; green, 150; blue, 190 }] (50,75) circle (1.5pt);
	\filldraw [color={rgb, 255:red, 37; green, 150; blue, 190 }] (100,75) circle (1.5pt);
	\filldraw [color={rgb, 255:red, 37; green, 150; blue, 190 }] (75,100) circle (1.5pt);
	
	\filldraw [color={rgb, 255:red, 37; green, 150; blue, 190 }] (50,125) circle (1.5pt);
	\filldraw [color={rgb, 255:red, 37; green, 150; blue, 190 }] (75,125) circle (1.5pt);
	\filldraw [color={rgb, 255:red, 37; green, 150; blue, 190 }] (100,125) circle (1.5pt);
	\filldraw [color={rgb, 255:red, 37; green, 150; blue, 190 }] (75,150) circle (1.5pt);
	
	\filldraw [color={rgb, 255:red, 37; green, 150; blue, 190 }] (125,100) circle (1.5pt);
	\filldraw [color={rgb, 255:red, 37; green, 150; blue, 190 }] (125,125) circle (1.5pt);
	\filldraw [color={rgb, 255:red, 37; green, 150; blue, 190 }] (150,125) circle (1.5pt);
	\filldraw [color={rgb, 255:red, 37; green, 150; blue, 190 }] (125,150) circle (1.5pt);

	\draw (115,15.0) node [anchor=north west][inner sep=0.75pt]  [font=\normalsize] [align=left] {$\displaystyle \mathbb{G}_w^j$};
	\draw (15,15) node [anchor=north west][inner sep=0.75pt]  [font=\normalsize,color={rgb, 255:red, 65; green, 117; blue, 5 }  ,opacity=1 ] [align=left] {$\displaystyle \{\tau_k,w_k\}_{k = 0}^{K}$};
	\draw (60,80) node [anchor=north west][inner sep=0.75pt]  [font=\normalsize,color={rgb, 255:red, 208; green, 2; blue, 27 }  ,opacity=1 ] [align=left] {$\displaystyle \mathbb{K}_w^j$};	
\end{tikzpicture}
}%
	\caption{Splitting process of 2D grid $\mathbb{G}_w^j$ (black) with additional points after update according to \cite{richter2023adaptive} (blue), due to the collection of waypoints in collision $\{\tau_k,w_k\}_{k = 0}^{K}$ (green) and those belonging to this grid $\mathbb{K}_w^j$ (red).}%
	\label{fig:splitting}%
\end{figure}
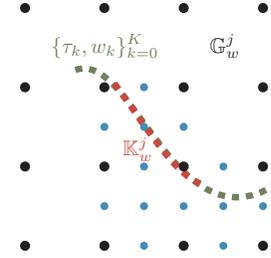

Let us focus now on the refinement evoked by an infeasible timestamped waypoint $\{\tau_k,w_k\}$.
First, we find the unique time index $j\in\N$ from $\mathbb{G}_\tau$ such that $\tau_{j-1} < \tau_k \leq \tau_j$.
Then, in order to maintain the Markov property for \eqref{eq:DPProb}, as recommended in \cite[Section IV.B]{britzelmeier2023dynamic}, both grids $j-1$ and $j$ will be refined, as in \GlobalPotato{} for the penalty adaptation.
Therefore, the set of cells selected for refinement is identified by collecting all cells containing at least an infeasible waypoint, namely
\begin{multline}
	\label{eq:SplittingSet}
	\mathbb{K}_w^j
	\coloneqq
	\left\{ c \,\vert\, \exists \{\tau,w\} \in \{\tau_k,w_k\}_{k = 0}^K, \right. \\ \left. c = C_w^j(w), \tau \in [\tau_{j-1},\tau_{j+1}] \right\} .
\end{multline}
Finally, the state grids for the next DP evaluation within \GlobalPotato{} are updated according to
$
	\mathbb{G}_w^j \gets \cup_{c \in \mathbb{K}_w^j} \mathcal{S}_w^j(c)
$.
While Figure~\ref{fig:splitting} illustrates the adaptive splitting process, Algorithm~\ref{alg:IterativeScheme} outlines the overall iterative procedure;
compare Figure~\ref{fig:flowchart} and \cite[Alg. 1]{britzelmeier2023dynamic}.

\begin{algorithm}[tbh]
	\SetAlgoLined
	\DontPrintSemicolon
	\caption{\GlobalPotato{} scheme with grid refinements}
	\label{alg:IterativeScheme}
	Initialize state and control grids $\{\mathbb{G}_w^j\}_{j=0}^M$, $\mathbb{G}_v$\;

\While{true}{
	
	$\{w_j, v_j,h\}_{j=0}^M \gets$ \eqref{eq:DPProb} with grids $\mathbb{G}_w^j$ \label{algstep:dp}\;
	
	\KwTry $\{x_j\}_{j=0}^M \gets \{w_j\}_{j=0}^M$ by \eqref{eq:SpaceMapping}\; 
	\quad{}build $\{p_j\}_{j=0}^M$\label{algstep:inverse_kinematics}\;
	
	\KwCatch $ \{w_k\}_{k = 0}^K \subset \{w_j\}_{j=0}^M$ with \eqref{eq:SpaceMapping} failed at $w_k$\;
	\quad{}\KwGoTo Step~\ref{algstep:grid_adaptation}\label{algstep:failed_inverse_kinematics}\;
	
	$\{t_i, x_i, u_i\}_{i=0}^N \gets$ \eqref{eq:LL_Problem} with $\{p_j\}_{j=0}^M$\label{algstep:nlp}\;
	
	\lIf {$\forall i \in \{1,\ldots,N\} ~g(x_i) \leq 0$}{\Return\label{algstep:constraintevaluation}}
	
	$\{w_k\}_{k = 0}^K \subset \{ \Omega(x_i) \}_{i=0}^N$ with $g(x_k) > 0$\;
	
	$\mathbb{G}_w^j \gets \cup_{c \in \mathbb{K}_w^j} \mathcal{S}_w^j(c)$ by \eqref{eq:SplittingSet}\label{algstep:grid_adaptation}\;
	
}
\end{algorithm}

\section{Numerical Tests}\label{sec:numerical}
In order to numerically demonstrate the capabilities of our approach, we apply it to the trajectory planning task of a robotic manipulator within a space application and compare the results to the original \GlobalPotato{} \cite{britzelmeier2023dynamic} as well as to a full discretization approach using Ipopt \cite{waechter2006implementation}.
The task is challenging due to the many degrees of freedom and the highly nonlinear dynamics.
Notice that a comparison with pure DP, including splitting or not, is impractical due to the curse of dimensionality.

\subsection{Implementation}
We implemented Algorithm~\ref{alg:IterativeScheme} and the original \GlobalPotato{} in \texttt{C++}.
Ipopt version 3.14.12 \cite{waechter2006implementation} is used, in combination with HSL linear solvers \cite{hsl} and the Ifopt interface version 2.1.3 \cite{ifopt}, to solve the auxiliary problems \eqref{eq:LL_Problem}, \eqref{eq:penalty_function}, and \eqref{algstep:inverse_kinematics}.
Further, we adopt the automatic differentiation tool \texttt{autodiff} v1.0.3 \cite{autodiff} to estimate the derivatives required when solving \eqref{eq:LL_Problem}.
To limit the cases where we observe nonzero penalties from local optima at feasible regions, we initially solve a problem without cost function holding solely the collision avoidance and mapping constraint and afterwards use its result as an initial guess on solving \eqref{eq:penalty_function}.
All the results presented were executed on an AMD Ryzen Threadripper 3990X system with 256~GB RAM.

\subsection{Problem}
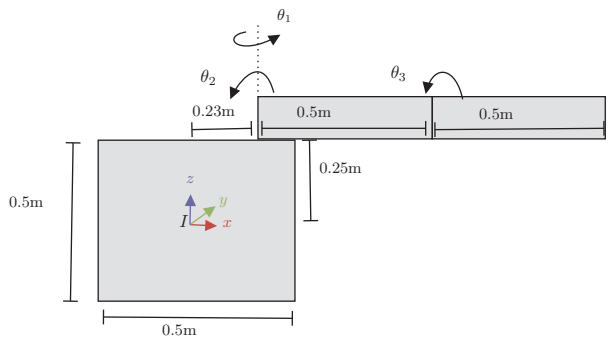
\begin{figure}[tb]
	\centering%
	\scalebox{0.666}{%

\tikzset{every picture/.style={line width=0.75pt}} 

\begin{tikzpicture}[x=0.75pt,y=0.75pt,yscale=-1,xscale=1]

\draw  [fill={rgb, 255:red, 155; green, 155; blue, 155 }  ,fill opacity=0.35 ] (100,129) -- (249,129) -- (249,251) -- (100,251) -- cycle ;
\draw [color={rgb, 255:red, 208; green, 2; blue, 27 }  ,draw opacity=1 ]   (169.5,193) -- (187,193.85) ;
\draw [shift={(190,194)}, rotate = 182.79] [fill={rgb, 255:red, 208; green, 2; blue, 27 }  ,fill opacity=1 ][line width=0.08]  [draw opacity=0] (8.93,-4.29) -- (0,0) -- (8.93,4.29) -- cycle    ;
\draw [color={rgb, 255:red, 126; green, 211; blue, 33 }  ,draw opacity=1 ]   (169.5,193) -- (186.56,180.75) ;
\draw [shift={(189,179)}, rotate = 144.32] [fill={rgb, 255:red, 126; green, 211; blue, 33 }  ,fill opacity=1 ][line width=0.08]  [draw opacity=0] (8.93,-4.29) -- (0,0) -- (8.93,4.29) -- cycle    ;
\draw [color={rgb, 255:red, 74; green, 83; blue, 226 }  ,draw opacity=1 ]   (169.5,193) -- (169.93,173) ;
\draw [shift={(170,170)}, rotate = 91.25] [fill={rgb, 255:red, 74; green, 83; blue, 226 }  ,fill opacity=1 ][line width=0.08]  [draw opacity=0] (8.93,-4.29) -- (0,0) -- (8.93,4.29) -- cycle    ;
\draw  [fill={rgb, 255:red, 155; green, 155; blue, 155 }  ,fill opacity=0.35 ] (221,96) -- (353,96) -- (353,128) -- (221,128) -- cycle ;
\draw  [fill={rgb, 255:red, 155; green, 155; blue, 155 }  ,fill opacity=0.35 ] (353,96) -- (484,96) -- (484,128) -- (353,128) -- cycle ;
\draw    (104,263) -- (247,264) ;
\draw [shift={(247,264)}, rotate = 180.4] [color={rgb, 255:red, 0; green, 0; blue, 0 }  ][line width=0.75]    (0,5.59) -- (0,-5.59)   ;
\draw [shift={(104,263)}, rotate = 180.4] [color={rgb, 255:red, 0; green, 0; blue, 0 }  ][line width=0.75]    (0,5.59) -- (0,-5.59)   ;
\draw    (81,131) -- (79,251) ;
\draw [shift={(79,251)}, rotate = 270.95] [color={rgb, 255:red, 0; green, 0; blue, 0 }  ][line width=0.75]    (0,5.59) -- (0,-5.59)   ;
\draw [shift={(81,131)}, rotate = 270.95] [color={rgb, 255:red, 0; green, 0; blue, 0 }  ][line width=0.75]    (0,5.59) -- (0,-5.59)   ;
\draw    (171,121) -- (216,120) ;
\draw [shift={(216,120)}, rotate = 178.73] [color={rgb, 255:red, 0; green, 0; blue, 0 }  ][line width=0.75]    (0,5.59) -- (0,-5.59)   ;
\draw [shift={(171,121)}, rotate = 178.73] [color={rgb, 255:red, 0; green, 0; blue, 0 }  ][line width=0.75]    (0,5.59) -- (0,-5.59)   ;
\draw    (260,129) -- (261,190) ;
\draw [shift={(261,190)}, rotate = 269.06] [color={rgb, 255:red, 0; green, 0; blue, 0 }  ][line width=0.75]    (0,5.59) -- (0,-5.59)   ;
\draw [shift={(260,129)}, rotate = 269.06] [color={rgb, 255:red, 0; green, 0; blue, 0 }  ][line width=0.75]    (0,5.59) -- (0,-5.59)   ;
\draw    (224,120) -- (348,120) ;
\draw [shift={(348,120)}, rotate = 180] [color={rgb, 255:red, 0; green, 0; blue, 0 }  ][line width=0.75]    (0,5.59) -- (0,-5.59)   ;
\draw [shift={(224,120)}, rotate = 180] [color={rgb, 255:red, 0; green, 0; blue, 0 }  ][line width=0.75]    (0,5.59) -- (0,-5.59)   ;
\draw    (355,121) -- (483,120) ;
\draw [shift={(483,120)}, rotate = 179.55] [color={rgb, 255:red, 0; green, 0; blue, 0 }  ][line width=0.75]    (0,5.59) -- (0,-5.59)   ;
\draw [shift={(355,121)}, rotate = 179.55] [color={rgb, 255:red, 0; green, 0; blue, 0 }  ][line width=0.75]    (0,5.59) -- (0,-5.59)   ;
\draw  [dash pattern={on 0.84pt off 2.51pt}]  (221,42) -- (221,96) ;
\draw    (212,47) .. controls (192.5,49.93) and (205.32,67.11) .. (235.63,50.36) ;
\draw [shift={(238,49)}, rotate = 149.3] [fill={rgb, 255:red, 0; green, 0; blue, 0 }  ][line width=0.08]  [draw opacity=0] (8.93,-4.29) -- (0,0) -- (8.93,4.29) -- cycle    ;
\draw    (201.37,95.1) .. controls (209.12,79.9) and (225.48,64.8) .. (233,93) ;
\draw [shift={(200,98)}, rotate = 293.63] [fill={rgb, 255:red, 0; green, 0; blue, 0 }  ][line width=0.08]  [draw opacity=0] (8.93,-4.29) -- (0,0) -- (8.93,4.29) -- cycle    ;
\draw    (348.75,90.92) .. controls (354.59,70.31) and (371.23,74.13) .. (376,98) ;
\draw [shift={(348,94)}, rotate = 281.77] [fill={rgb, 255:red, 0; green, 0; blue, 0 }  ][line width=0.08]  [draw opacity=0] (8.93,-4.29) -- (0,0) -- (8.93,4.29) -- cycle    ;

\draw (193,188.4) node [anchor=north west][inner sep=0.75pt]  [font=\small]  {$\textcolor[rgb]{0.82,0.01,0.11}{x}$};
\draw (190,170.4) node [anchor=north west][inner sep=0.75pt]  [font=\small]  {$\textcolor[rgb]{0.49,0.83,0.13}{y}$};
\draw (165,154.4) node [anchor=north west][inner sep=0.75pt]  [font=\small]  {$\textcolor[rgb]{0.29,0.31,0.89}{z}$};
\draw (147,267.4) node [anchor=north west][inner sep=0.75pt]  [font=\small]  {$0.5\si{m}$};
\draw (31,170.4) node [anchor=north west][inner sep=0.75pt]  [font=\small]  {$0.5\si{m}$};
\draw (170,101.4) node [anchor=north west][inner sep=0.75pt]  [font=\small]  {$0.23\si{m}$};
\draw (266,144.4) node [anchor=north west][inner sep=0.75pt]  [font=\small]  {$0.25\si{m}$};
\draw (249,102.4) node [anchor=north west][inner sep=0.75pt]  [font=\small]  {$0.5\si{m}$};
\draw (387,103.4) node [anchor=north west][inner sep=0.75pt]  [font=\small]  {$0.5\si{m}$};
\draw (160,184.4) node [anchor=north west][inner sep=0.75pt]    {$I$};
\draw (234,28.4) node [anchor=north west][inner sep=0.75pt]  [font=\small]  {$\theta _{1}$};
\draw (176,74.4) node [anchor=north west][inner sep=0.75pt]    {$\theta _{2}$};
\draw (320,71.4) node [anchor=north west][inner sep=0.75pt]    {$\theta _{3}$};

\end{tikzpicture}%
}%
%
	\caption{Illustration of the satellite body with robotic arm. The coordinate axes of the satellites inertial frame $I$ are highlighted.}%
	\label{fig:satellite}%
\end{figure}

We consider the model of a free-floating satellite with a two link robotic arm attached, see Figure~\ref{fig:satellite}.
Its $9$-dimensional configuration $\tilde{q} = \{q,\theta\}$ comprises the satellite body's pose $q \in \R^6$, with respect to the inertial reference frame $I$, and the joint angles of the robotic arm $\theta \in \R^3$.
The satellite-arm system is acted upon with controls $\tilde{u} \in \R^9$, each one associated to a degree of freedom.
Instead of considering a fully free-floating system, we simplify the robotic arm dynamics, approximating them with double-integrators, as $\ddot{\theta}_i = u_{i+6}$ for $i = 1,2,3$.
Each body's pose with respect to $I$ is determined using the robotic arm's forward kinematics.
Subsequently, (translational and rotational) velocities and accelerations are obtained using the first and second time derivatives.
For the remaining degrees of freedom, Lagrange mechanics yield the equations of motion \cite{bertolazzi2007realtime} in the form
\begin{equation*}
	M \left(\tilde{q},\dot{\tilde{q}}\right) \ddot{q}
	=
	\tilde{f} \left(\tilde{q},\dot{\tilde{q}},\tilde{u}_q\right)
	,
\end{equation*}
with mass matrix $M \in \R^{6 \times 6}$ and satellite control signals $\tilde{u}_q \coloneqq \{u_1,\ldots,u_6\}$.

In view of \eqref{eq:OverallProb}, we define the state vector of our second-order satellite model as $x \coloneqq \{\tilde{q},\dot{\tilde{q}}\} \in \R^{18}$, with control $u \coloneqq \tilde{u} \in \R^9$ and include the equations of motion into the multi dimensional dynamics function, defined component wise as
\begin{align}
	&f_i(x(t),u(t)) = x(t)_{i+9} ~ && i = 1,\ldots,9 \notag \\
	&f_i(x(t),u(t))) = M^{-1}\tilde{f}(x(t),u(t)) \notag ~&& i = 10,\ldots,15 \\
	&f_i(x(t),u(t)) = u(t)_{i-9} ~&& i = 16,17,18 \notag	 
\end{align}  
To complete the geometric information provided in Figure~\ref{fig:satellite} we fix the satellite body's depth to $0.4 \si{m}$, its mass to $50 \si{kg}$ and its inertia diagonal matrix values to $\{1.71, 1.71, 2.08\} \si{kg m^2}$.
Accordingly, the  cylindric robotic links have a radius of $0.05\si{m}$, diagonal inertia entries $\{0.006, 0.107, 0.107\} \si{kg m^2}$ and a mass of $5\si{kg}$ each.
We limit the state and control according to Table \ref{tab:Limits} and the time bounds with $T_{\min} = 0 \si{s}$ and $T_{\max} = 500\si{s}$.

\begin{table}[tbh]
	\centering%
	\caption{State and Control boundaries.}%
	\label{tab:Limits}%
	\setlength{\tabcolsep}{2pt}
	\begin{tabular}{cccccccccc}
		\toprule
		$\tilde{q}_{\min}$ & $\tilde{q}_{\max}$ & $\tilde{q}'_{\min}$ &  $\tilde{q}'_{\max}$ & $\tilde{u}_{\min}$ & $\tilde{u}_{\max}$ & $w_{min}$ & $w_{max}$ & $v_{min}$ & $v_{max}$\\
		\midrule
	 $[\si{m}]$ & $[\si{m}]$ & $[\si{m}/\si{s}]$ & $[\si{m}/\si{s}]$ &  $[\si{N}]$ &  $[\si{N}]$&$[\si{m}]$ & $[\si{m}]$ &$[\si{m}/\si{s}]$&$[\si{m}/\si{s}]$\\
	 -1.0 & 15.0 & -0.5 & 0.5 & -1.2 & 1.2 & -1.0 & 15.0 & -0.04 & 0.04\\
		-5.0 & 5.0 & -0.5 & 0.5 & -1.2 & 1.2 & -5.0 & 5.0 & -0.04 & 0.04\\
		-1.0 & 15.0 &-0.5 & 0.5 & -1.2 & 1.2 & -1.0 & 15.0 & -0.04 & 0.04\\
		$[\si{rad}]$ & $[\si{rad}]$ & $[\si{rad}/\si{s}]$ & $[\si{rad}/\si{s}]$ &  $[\si{Nm}]$ &  $[\si{Nm}]$&&&&\\
		-$\infty$ &$\infty$&-$\pi/6$ & $\pi/6$ & -1.5 & 1.5 &&&&\\
		$-1.47$& $1.47$ &-$\pi/6$& $\pi/6$  & -1.5 & 1.5 &&&&\\
		-$\infty$ & $\infty$ &-$\pi/6$ & $\pi/6$ & -1.5 & 1.5  &&&&\\
		-$\pi/2$ & $\pi/2$ & -$\pi/10$ & $\pi/10$ & -1.0 & 1.0 &&&& \\
		-$\pi$ & 0.0 &-$\pi/10$ & $\pi/10$ & -1.0  & 1.0  &&&&\\
		-$\pi/2$ & $\pi/2$& -$\pi/10$& $\pi/10$ & -1.0 & 1.0 &&&&\\
		\bottomrule
	\end{tabular}%
\end{table}

In addition, we consider four rectangular static obstacles in our environment, illustrated in Figure \ref{fig:SatelliteTrajectories}.
We formulate the collision avoidance constraint as in \cite[Eq. 13]{britzelmeier2023dynamic}, based on the signed distance between polygons and using a safety bound $\varepsilon = 0.01 \si{m}$.
We formulate the final conditions as
\begin{align*}
	&b_0(x) = \alpha \max\{ \| \Omega(x) - w_{\text{ref}} \| - r, 0 \} \\
	&b_i(x) = x_{i+9} && i = 1,\ldots,9
\end{align*}
with $b_0$ returning the positive distance to the goal TCP point $w_{\text{ref}} = \{12.0,0.0,1.0\}\si{m}$, scaled by parameter $\alpha = 10^3$, and zero on a ball with radius $r = 1 \si{m}$ around it; see Figure~\ref{fig:SatelliteTrajectories}.
The other components of $b$ demand zero final velocity.

\subsection{Setup}
To solve the above problem using Algorithm~\ref{alg:IterativeScheme} we choose function $\Omega$ to represent the forward kinematics of the satellite, mapping a 18-dimensional configuration to the corresponding 3-dimensional tool center point (TCP) position.
We initially discretize the DP state and control space within the bounds presented in Table~\ref{tab:Limits}
into $10$ points along each dimension as well as 10 points on the discretization of the time step size $h_i$.
We set the dynamics function to $\bar{f}(w,v) = v$, the penalty scaling factor $\rho = 40$ and consider $M = 20$ DP time steps.
For the inverse space mapping \eqref{eq:SpaceMapping}, we set parameters $\varrho_1 = \varrho_2 = 1$.
We follow the arising waypoints using $p_j(x) = \Omega(x) - w_j$ as well as a minimum effort objective $J(x,u,T) = \int_{0}^{T_{\max}} \|u(t)\|_2^2 \mathrm{d}t$ in \eqref{eq:LL_Problem}, discretizing with an equidistant grid with $N = 200$ time intervals, piecewise constant control parametrization, and finite difference approximation of the dynamics with trapezoidal rule.
To realize the adaptive grid refinement, we use the cube-based grid structure proposed in \cite{richter2023adaptive} along with its respective cell identification and splitting operators.
In the experiments, we took the lower part of the initial state space grid e.g. $\{10,10,4\}$ points and tried to solve \eqref{eq:SpaceMapping} on each of them.
On success, we used the resulting configuration $\tilde{q}_0$ to form an initial condition $x_0 \coloneqq\{\tilde{q}_0,0\}$ and solved the respective problem using (at most) 15 iterations of Algorithm~\ref{alg:IterativeScheme} and \GlobalPotato{}.
After adopting the discrete-time version of \eqref{eq:DPProb} without the stepsize in the objective approximation, as in \cite[Eq. 6]{britzelmeier2023dynamic}, we repeated these experiments with \eqref{eq:DiscDPProb} for both Algorithm~\ref{alg:IterativeScheme} and \GlobalPotato{}.

Finally, to compare against a full discretization approach, we take 6 samples, distributed among the space of initial positions, where Algorithm~\ref{alg:IterativeScheme} was able to terminate (on the first set of experiments).
Discretizing \eqref{eq:OverallProb} as for \eqref{eq:LL_Problem}, the resulting zero-objective NLP is given to Ipopt,
initialized either from a zero guess or with a collision-free path obtained using the RRTConnect implementation from OMPL \cite{sucan2012the-open-motion-planning-library}.

\newenvironment{customlegend}[1][]{%
	\begingroup
	\csname pgfplots@init@cleared@structures\endcsname
	\pgfplotsset{#1}%
}{%
	\csname pgfplots@createlegend\endcsname
	\endgroup
}%
\def\addlegendimage{\csname pgfplots@addlegendimage\endcsname}

\begin{figure}[tbh]
	\centering
	\begin{tikzpicture}
		\begin{customlegend}[legend columns=2,legend style={font=\footnotesize,align=left,draw=none,column sep=2ex},
			legend entries={Initial Grid, Waypoints Iter 0, Added Grid Points Iter 0, Waypoints Iter 1, Added Grid Points Iter 1, Waypoints Iter 2}]
			\addlegendimage{only marks, mark=+,black}
			\addlegendimage{only marks, mark=*,mydarkblue}
			\addlegendimage{only marks, mark=+,mydarkblue}
			\addlegendimage{only marks, mark=*,mylightgreen}
			\addlegendimage{only marks, mark=+,mylightgreen}
			\addlegendimage{only marks, mark=*,mydarkgreen}
		\end{customlegend}
	\end{tikzpicture}
	\includegraphics[width = 0.75\linewidth]{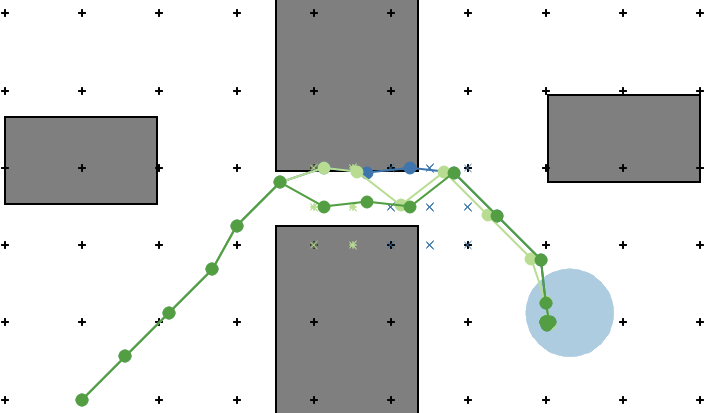}
	\caption{Illustration of state space grids and waypoints at DP iteration 6 resulting from iterations 1, 2, and 3 of Algorithm~\ref{alg:IterativeScheme} using initial condition $\tilde{q}_0 = \{-0.2,-0.14,-0.68,-0.49,0.34,0.43,0.13,-0.02,0.1\}$. The goal region (light blue) as well as the obstacles (gray) are highlighted.}
	\label{fig:SatelliteTrajectories}	
\end{figure}

\subsection{Results and Discussion}
The waypoints resulting from the first three iterations of Algorithm~\ref{alg:IterativeScheme} and a certain initial condition
are illustrated in Figure~\ref{fig:SatelliteTrajectories}.
In the first iteration, since there are no collision-free grid points in the narrow passage, the points cross the obstacles at the region of smallest penalty, namely where there is minimum overlap.
As shown in Figure~\ref{fig:SatelliteTrajectories} the addition of state space points quickly results in a sampling of the gap between the obstacles, leading to feasible waypoints at iteration 2 and finally to a collision-free trajectory.
In contrast, the original \GlobalPotato{} was not able to provide a feasible solution.
Evaluating 254 initial conditions in total, \GlobalPotato{} provided solutions in 132, and our algorithm in 163 cases. 
We did not observe any case where \GlobalPotato{} was able to find a solution but our algorithm was not.
Similar (albeit not identical) results were obtained in the experiments with the discrete-time version \eqref{eq:DiscDPProb} of \eqref{eq:DPProb}.
Finally, Ipopt was able to solve only 4 of the 6 problems when starting from a zero initial guess, while it managed to solve all of them when initialized with collision-free paths.
We therefore numerically demonstrated the improved efficiency and robustness of our algorithm, due to its enhanced capability to provide feasible solutions, compared to the original \GlobalPotato{} and a trivially initialized full discretization approach.

\section{Conclusions}\label{sec:conclusions}
In this work, we proposed a method combining dynamic programming, direct optimal control techniques,
and adaptive refinement of state space discretizations to solve a trajectory planning problem with obstacle avoidance.
Compared to adopting these techniques separately, enhanced robustness and efficiency of our method stem from the grid refinements, leading to improved convergence and computational scalability.
These achievements are numerically demonstrated on a trajectory planning task for a satellite manipulator system. 

\bibliographystyle{IEEEtran}
\bibliography{IEEEabrv,references}

\begin{thebibliography}{10}
\providecommand{\url}[1]{#1}
\csname url@rmstyle\endcsname
\providecommand{\newblock}{\relax}
\providecommand{\bibinfo}[2]{#2}
\providecommand\BIBentrySTDinterwordspacing{\spaceskip=0pt\relax}
\providecommand\BIBentryALTinterwordstretchfactor{4}
\providecommand\BIBentryALTinterwordspacing{\spaceskip=\fontdimen2\font plus
\BIBentryALTinterwordstretchfactor\fontdimen3\font minus
  \fontdimen4\font\relax}
\providecommand\BIBforeignlanguage[2]{{%
\expandafter\ifx\csname l@#1\endcsname\relax
\typeout{** WARNING: IEEEtran.bst: No hyphenation pattern has been}%
\typeout{** loaded for the language `#1'. Using the pattern for}%
\typeout{** the default language instead.}%
\else
\language=\csname l@#1\endcsname
\fi
#2}}

\bibitem{hoy2015algorithms}
M.~Hoy, A.~S. Matveev, and A.~V. Savkin, ``Algorithms for collision-free
  navigation of mobile robots in complex cluttered environments: a survey,''
  \emph{Robotica}, vol.~33, no.~3, pp. 463--497, 2015.

\bibitem{bellman1957dynamic}
R.~Bellman, \emph{Dynamic Programming}.\hskip 1em plus 0.5em minus 0.4em\relax
  Princeton, New Jersey: Princeton University Press, 1957.

\bibitem{bertsekas2005dynamic}
D.~P. Bertsekas, \emph{Dynamic Programming \& Optimal Control}.\hskip 1em plus
  0.5em minus 0.4em\relax Belmont, Massachusetts: Athena Scientific, 2005, 3rd
  edition.

\bibitem{cervellera2013quasi}
C.~Cervellera, M.~Gaggero, D.~Macci{\`o}, and R.~Marcialis, ``Quasi-random
  sampling for approximate dynamic programming,'' in \emph{The 2013
  International Joint Conference on Neural Networks (IJCNN)}, 2013, pp. 1--8.

\bibitem{weber2017optimized}
A.~Weber, M.~Rungger, and G.~Reissig, ``Optimized state space grids for
  abstractions,'' \emph{IEEE Transactions on Automatic Control}, vol.~62,
  no.~11, pp. 5816--5821, 2017.

\bibitem{richter2023adaptive}
R.~Richter, A.~Britzelmeier, and M.~Gerdts, ``An adaptive mesh dynamic
  programming algorithm for robotic manipulator trajectory planning,'' in
  \emph{2023 European Control Conference (ECC)}, Bucharest, Romania, 2023, pp.
  1--8.

\bibitem{munos2002variable}
R.~Munos and A.~Moore, ``Variable resolution discretization in optimal
  control,'' \emph{Machine Learning}, vol.~49, no.~2, pp. 291--323, 2002.

\bibitem{elbanhawi2014sampling}
M.~Elbanhawi and M.~Simic, ``Sampling-based robot motion planning: A review,''
  \emph{IEEE access}, vol.~2, pp. 56--77, 2014.

\bibitem{lavalle1998rapidly}
S.~M. LaValle, ``Rapidly-exploring random trees: A new tool for path
  planning,'' Department of Computer Science, Iowa State University, Tech. Rep.
  Research Report 9811, 1998.

\bibitem{gerdts2023optimal}
M.~Gerdts, \emph{Optimal Control of {ODEs} and {DAEs}}.\hskip 1em plus 0.5em
  minus 0.4em\relax De Gruyter Oldenbourg, 2023, 2nd Edition.

\bibitem{waechter2006implementation}
A.~W{\"a}chter and L.~T. Biegler, ``On the implementation of an interior-point
  filter line-search algorithm for large-scale nonlinear programming,''
  \emph{Mathematical Programming}, vol. 106, no.~1, pp. 25--57, 2006.

\bibitem{bueskens2013esa}
C.~B{\"u}skens and D.~Wassel, ``The {ESA} {NLP} solver {WORHP},'' in
  \emph{Modeling and Optimization in Space Engineering}, G.~Fasano and J.~D.
  Pintér, Eds.\hskip 1em plus 0.5em minus 0.4em\relax Springer New York, 2013,
  vol.~73, pp. 85--110.

\bibitem{zhang2021optimization}
X.~Zhang, A.~Liniger, and F.~Borrelli, ``Optimization-based collision
  avoidance,'' \emph{IEEE Transactions on Control Systems Technology}, vol.~29,
  no.~3, pp. 972--983, 2021.

\bibitem{guthrie2022differentiable}
J.~Guthrie, ``A differentiable signed distance representation for continuous
  collision avoidance in optimization-based motion planning,'' in \emph{2022
  IEEE 61st Conference on Decision and Control (CDC)}, 2022, pp. 7214--7221.

\bibitem{zucker2013chomp}
M.~Zucker, N.~Ratliff, A.~D. Dragan, M.~Pivtoraiko, M.~Klingensmith, C.~M.
  Dellin, J.~A. Bagnell, and S.~S. Srinivasa, ``{CHOMP}: Covariant
  {H}amiltonian optimization for motion planning,'' \emph{The International
  Journal of Robotics Research}, vol.~32, no. 9--10, pp. 1164--1193, 2013.

\bibitem{schulman2013finding}
J.~Schulman, J.~Ho, A.~Lee, I.~Awwal, H.~Bradlow, and P.~Abbeel, ``Finding
  locally optimal, collision-free trajectories with sequential convex
  optimization,'' in \emph{Proceedings of Robotics: Science and Systems},
  Berlin, Germany, 2013.

\bibitem{natarajan2021inter}
R.~Natarajan, H.~Choset, and M.~Likhachev, ``Interleaving graph search and
  trajectory optimization for aggressive quadrotor flight,'' \emph{IEEE
  Robotics and Automation Letters}, vol.~6, no.~3, pp. 5357--5364, 2021.

\bibitem{britzelmeier2023dynamic}
A.~Britzelmeier, A.~De~Marchi, and R.~Richter, ``Dynamic and nonlinear
  programming for trajectory planning,'' \emph{IEEE Control Systems Letters},
  vol.~7, pp. 2569--2574, 2023.

\bibitem{marcucci2022motion}
T.~Marcucci, M.~Petersen, D.~von Wrangel, and R.~Tedrake, ``Motion planning
  around obstacles with convex optimization,'' \emph{Science Robotics}, vol.~8,
  no.~84, p. eadf7843, 2023.

\bibitem{britzelmeier2021decomposition}
A.~Britzelmeier and A.~Dreves, ``A decomposition algorithm for {N}ash
  equilibria in intersection management,'' \emph{Optimization}, vol.~70,
  no.~11, pp. 2441--2478, 2021.

\bibitem{hsl}
\BIBentryALTinterwordspacing
HSL, ``A collection of fortran codes for large scale scientific computation,''
  2013. [Online]. Available: \url{http://www.hsl.rl.ac.uk}
\BIBentrySTDinterwordspacing

\bibitem{ifopt}
\BIBentryALTinterwordspacing
A.~W. Winkler, ``{Ifopt} - a modern, light-weight, {Eigen}-based {C++}
  interface to nonlinear programming solvers {IPOPT} and {SNOPT}.'' 2018.
  [Online]. Available: \url{https://github.com/ethz-adrl/ifopt}
\BIBentrySTDinterwordspacing

\bibitem{autodiff}
\BIBentryALTinterwordspacing
A.~M.~M. Leal, ``autodiff, a modern, fast and expressive {C++} library for
  automatic differentiation,'' 2018. [Online]. Available:
  \url{https://autodiff.github.io}
\BIBentrySTDinterwordspacing

\bibitem{bertolazzi2007realtime}
E.~Bertolazzi, F.~Biral, and M.~Da~Lio, ``Real-time motion planning for
  multibody systems: Real life application examples,'' \emph{Multibody System
  Dynamics}, vol.~17, no. 2--3, pp. 119--139, 2007.

\bibitem{sucan2012the-open-motion-planning-library}
I.~A. {\c{S}}ucan, M.~Moll, and L.~E. Kavraki, ``The {O}pen {M}otion {P}lanning
  {L}ibrary,'' \emph{{IEEE} Robotics \& Automation Magazine}, vol.~19, no.~4,
  pp. 72--82, December 2012, \url{https://ompl.kavrakilab.org}.

\end{thebibliography}

\end{document}